\documentclass[11pt,a4paper,twoside,reqno]{amsart}

\usepackage[a4paper,margin=1in]{geometry}
\usepackage{amsmath,amssymb,amsthm,mathtools}
\usepackage{microtype}
\usepackage[numbers,sort&compress]{natbib}
\usepackage[hidelinks]{hyperref}

\numberwithin{equation}{section}

\newtheorem{theorem}{Theorem}[section]
\newtheorem{proposition}[theorem]{Proposition}
\newtheorem{lemma}[theorem]{Lemma}
\newtheorem{corollary}[theorem]{Corollary}

\theoremstyle{definition}

\theoremstyle{remark}

\newtheorem{openproblem}[theorem]{Open Problem}

\newcommand{\E}{\mathbb E}
\newcommand{\PP}{\mathbb P}
\newcommand{\1}{\mathbf 1}
\newcommand{\normT}[1]{\lVert #1\rVert_{\mathbb T}}
\newcommand{\dd}{\,\mathrm d}

\hypersetup{
  pdftitle={Poissonian Pair Correlation and Discrepancy Rates},
  pdfauthor={},
  pdfsubject={Poissonian pair correlation and one-dimensional star discrepancy},
  pdfkeywords={Poissonian pair correlation, star discrepancy, low-discrepancy sequence}
}

\title{Poissonian Pair Correlation and Discrepancy Rates}
\author{Bo Tan, \ \ Zhen Zhang$^{\ast}$, \ \ Qing-Long Zhou}
\thanks{$^{\ast}$Corresponding author
}
\date{}

\begin{document}

\begin{abstract}
We prove the existence of a sequence in $[0,1)$ with Poissonian pair
correlation (PPC) and optimal-order star discrepancy
$D_N^*=O((\log N)/N)$, and provide a deterministic construction
by applying a conditional-expectation rule to a nested dyadic random
model. We also show that PPC imposes no universal discrepancy decay
rate: for every positive sequence $\rho_N\to0$, there exists a PPC
sequence with $D_N^*/\rho_N\to\infty$. Under suitable regularity and
growth assumptions, we realize prescribed discrepancy asymptotics,
including $D_N^*\sim c_\beta N^{-\beta}$ for every $0<\beta<1$, with
$c_\beta>0$. We further construct a PPC sequence satisfying
$D_N^*\asymp(\log N)/N$ for all sufficiently large $N$.


\end{abstract}

\maketitle

\section{Introduction}

The theory of uniform distribution modulo one was developed in the
work of Weyl \cite{Weyl1916}. A sequence
$(x_n)_{n\ge1}\subset[0,1)$ is uniformly distributed modulo one if
\[
  \frac1N\sharp\{1\le n\le N:x_n\in[a,b)\}
  \to  b-a
  \qquad(N\to\infty)
\]
for every $0\le a<b\le1$. A basic example is the Kronecker sequence
$(\{n\alpha\})_{n\ge1}$ for any irrational $\alpha$, where
$\{t\}=t-\lfloor t\rfloor$ denotes the fractional part of $t$.

For  $x=(x_n)\subset[0,1)$, put
\[
A_N^x(t)=\sharp\{1\le n\le N:x_n<t\}
\qquad(N\ge1,\ 0\le t\le1),
\]
with $A_0^x(t)=0$.
The star discrepancy of the first $N$ points is
\begin{equation}
\notag
  D_N^*= D_N^*(x):=\sup_{0\le t\le1}
  \left|\frac{A_N^x(t)}N-t\right|.
\end{equation}
A sequence $x$ is uniformly distributed modulo one if and only if
$D_N^*(x)\to0$ as $N\to\infty$; see
\citet[Chapter~2, Corollary~1.1]{KuipersNiederreiter1974}.
Thus star discrepancy measures the rate of convergence to uniform
distribution. A sequence is called a low-discrepancy sequence if
\(D_N^*=O((\log N)/N)\). This order is optimal for infinite
one-dimensional sequences: Schmidt \cite{Schmidt1972} proved that
\(ND_N^*\gg\log N\) for infinitely many \(N\) for every infinite sequence in
\([0,1)\). Both the Kronecker sequence \((\{n\alpha\})_{n\ge1}\) with
\(\alpha\) having bounded partial quotients in its continued fraction
expansion, and the van der Corput sequence, attain this order; see
\citet[Chapter~2, Theorems~3.4 \& 3.5]{KuipersNiederreiter1974}. We also refer the
reader to this book for further background and history, as well as for more
information.

Pair correlation concerns distances between pairs of points at the scale
\(1/N\). We identify \([0,1)\) with \(\mathbb T=\mathbb R/\mathbb Z\) and write
\(\|\cdot\|_{\mathbb{T}}=\min_{m\in\mathbb Z}|\cdot-m|\) for the distance to the nearest
integer. For a sequence \(x=(x_n)\), we define
\begin{equation}
\notag
 R_N(s)= R_N^x(s):=\frac1N
  \sharp\left\{1\le i\ne j\le N:
  \normT{x_i-x_j}\le\frac{s}{N}\right\}
  \qquad (s>0).
\end{equation}
The sequence has Poissonian pair correlation (PPC) if, for every \(s>0\),
\begin{equation}
\notag
  R_N(s)\to  2s
  \qquad(N\to\infty).
\end{equation}
 This is the almost-sure limit for independent uniform
random points on \(\mathbb T\).

El-Baz, Marklof and Vinogradov \cite{ElBazMarklofVinogradov2015}
proved that $(\{\sqrt n\})_{n\ge1}$ has PPC after the terms with
square $n$ are omitted. Metric examples include
$(\{n^2\alpha\})_{n\ge1}$ and $(\{2^n\alpha\})_{n\ge1}$,
which have PPC for Lebesgue-almost every $\alpha$, by results of
Rudnick and Sarnak \cite{RudnickSarnak1998} and
Rudnick and Zaharescu \cite{RudnickZaharescu1999}, respectively.

Uniform distribution concerns frequencies in fixed intervals, whereas
PPC concerns pair counts at the shrinking scale $s/N$.
PPC implies uniform distribution, as proved independently by
Aistleitner, Lachmann and Pausinger
\cite{AistleitnerLachmannPausinger2018} and by Grepstad and Larcher
\cite{GrepstadLarcher2017}.
The converse fails: for irrational $\alpha$, the Kronecker sequence
$(\{n\alpha\})_{n\ge1}$ is uniformly distributed but does not have PPC.
Even low discrepancy does not imply PPC, as the van der Corput
sequence shows; both examples are covered by the results of
\citet{LarcherStockinger2020}.

The relation between PPC and discrepancy also has a quantitative
form. Grepstad and Larcher \cite{GrepstadLarcher2017} and
Steinerberger \cite{Steinerberger2018} obtained discrepancy bounds
from quantitative estimates for pair-correlation statistics.

We study which discrepancy rates are compatible with PPC.
We first construct a sequence with PPC and optimal-order star
discrepancy\footnote{While writing up the present manuscript, we became aware
of related work on arXiv by H.~Porath,
\emph{On low-discrepancy sequences and Poissonian pair correlation}
(\href{https://arxiv.org/abs/2609.21689}{arXiv:2609.21689}),
which also established the existence of low-discrepancy sequences
with Poissonian pair correlation. We note that our construction
of the sequence is fundamentally distinct from that in the
aforementioned work.}. The construction yields a random sequence with low discrepancy for every
realization and PPC almost surely, together with a deterministic
conditional-expectation rule that selects one such realization.
We also obtain a PPC sequence satisfying
$D_N^*\asymp(\log N)/N$ for all sufficiently large $N$.
More generally, we construct PPC sequences with prescribed discrepancy
asymptotics $D_N^*\sim c_\psi\psi(N)$ under suitable regularity and
growth assumptions on $\psi$, where $c_\psi>0$.
These include every polynomial rate $N^{-\beta}$ with $0<\beta<1$.
We also show that, for every positive sequence $\rho_N\to0$, there
exists a PPC sequence with $D_N^*/\rho_N\to\infty$.
Thus PPC imposes no universal discrepancy decay rate.

 Throughout the paper, we use Landau and Vinogradov notation interchangeably. For instance,  we write \(a_N\ll b_N\), or equivalently \(a_N=O(b_N)\), if there exists a constant \(C>0\) such that \(b_N>0\) and \(|a_N|\le C b_N\) for all sufficiently large \(N\); sometimes a subscript on \(O\) or \(\ll\) is added to indicate parameters on which the implied constant \(C\) may depend. We write \(a_N\asymp b_N\) if both \(a_N\ll b_N\) and \(b_N\ll a_N\), and we write \(a_N\sim b_N\) if \(a_N/b_N\to1\) as \(N\to\infty\). We write \(a_N=o(b_N)\), or equivalently \(b_N=\omega(a_N)\), if \(a_N/b_N\to0\) as \(N\to\infty\).

We denote the indicator of a set $A$ by $\1_A$.
The notation $\mathcal B(S)$ denotes the Borel $\sigma$-algebra of $S$.
The symbols $\PP$ and $\E$ denote probability and expectation
for the random construction under consideration.
Let $\lambda$ denote normalized Lebesgue measure on $\mathbb T$.
All integrals on $\mathbb T$ are taken with respect to this measure.


Our main results read as follows.


\begin{theorem}
\label{thm:low-discrepancy-background}
There exists a deterministic sequence $(y_n)_{n\ge1}\subset[0,1)$ such that
\[
 R_N^y(s)\to  2s
 \qquad(N\to\infty)
 \quad\text{for every }s>0
\]
and
\begin{equation}
\label{eq:main-low-discrepancy}
 ND_N^*(y_1,\ldots,y_N)
 \le\frac{5}{2\log2}\log N
 \qquad(N\ge2).
\end{equation}
A deterministic conditional-expectation rule selects such a sequence
as a realization of a random sequence satisfying both properties
almost surely.
\end{theorem}

\begin{theorem}
\label{thm:critical-order}
There exists a sequence
\(
(x_n)_{n\ge1}\subset\mathbb T
\)
with PPC such that
\begin{equation}
\notag
D_N^*\asymp \frac{\log N}{N}
\end{equation}
for all sufficiently large \(N\).
\end{theorem}


\begin{theorem}
\label{thm:prescribed-intermediate}
Let $\psi:\mathbb N\to(0,\infty)$, and put
\(
h_N=N\psi(N)\) with \(h_0=0.\)  
Assume that, for some $L>0$,
\begin{equation}
\label{eq:intermediate-increments}
0\le h_N-h_{N-1}\le L
\qquad(N\ge1),
\end{equation}
\begin{equation}
\label{eq:intermediate-growth-lower}
h_N=\omega(\log N) 
\qquad(N\to\infty)
\end{equation}
and
\begin{equation}
\label{eq:intermediate-growth-upper}
h_N\ll\sqrt{N\log N}.
\end{equation}
Then there exist a sequence
$(x_n)_{n\ge1}\subset\mathbb T$ and a constant $c_\psi>0$ such that:
\begin{enumerate}
\item the sequence $(x_n)$ has PPC;

\item its star discrepancy satisfies
\begin{equation}
\notag
ND_N^*
=
c_\psi h_N
+
O\left(\log N+\sqrt{h_N\log h_N}\right)
\end{equation}
for all sufficiently large $N$. In particular,
\begin{equation}
\notag
D_N^*\sim c_\psi\psi(N)
\qquad(N\to\infty).
\end{equation}
\end{enumerate}
\end{theorem}

\begin{theorem}
\label{thm:prescribed-above-square-root}
Let $\psi:\mathbb N\to(0,\infty)$, and put
\(
h_N=N\psi(N)\) with \(h_0=0.\)
Assume that, for some $L>0$,
\begin{equation}
\notag
0\le h_N-h_{N-1}\le L
\qquad(N\ge1),
\end{equation}
and that
\begin{equation}
\label{eq:above-empirical-scale}
h_N=\omega(\sqrt{N\log N})
\qquad(N\to\infty)
\end{equation}
and
\begin{equation}
\notag
h_N=o(N)
\qquad(N\to\infty).
\end{equation}
Then there exist a sequence
$(x_n)_{n\ge1}\subset\mathbb T$ and a constant $c_\psi>0$ such that:
\begin{enumerate}
\item the sequence $(x_n)$ has PPC;

\item its star discrepancy satisfies
\begin{equation}
\notag
ND_N^*
=c_\psi h_N
+O\left(\sqrt{N\log N}\right)
\end{equation}
for all sufficiently large $N$. In particular,
\begin{equation}
\notag
D_N^*\sim c_\psi\psi(N)
\qquad(N\to\infty).
\end{equation}
\end{enumerate}
\end{theorem}

We record several immediate consequences of the main theorems, with proofs
that assume the theorems.

\begin{corollary}
For every $0<\beta<1$, there exist a PPC sequence and a constant
$c_\beta>0$ such that $D_N^*\sim c_\beta N^{-\beta}$ as $N\to\infty$.
\end{corollary}

\begin{proof}
Apply Theorem~\ref{thm:prescribed-intermediate} for $1/2\le\beta<1$
and Theorem~\ref{thm:prescribed-above-square-root} for $0<\beta<1/2$,
in both cases with $\psi(N)=N^{-\beta}$.
\end{proof}

\begin{corollary}
For every positive sequence $(\rho_N)_{N\ge1}$ with $\rho_N\to0$,
there exists a PPC sequence such that
\[
\frac{D_N^*}{\rho_N}\to\infty
\qquad(N\to\infty).
\]
Thus there is no positive sequence $\rho_N\to0$ such that
$D_N^*=O(\rho_N)$ for every PPC sequence.
\end{corollary}

\begin{proof}
Put
\[
u_n=\sup_{m\ge n}\sqrt{\rho_m}+n^{-1/4},
\qquad
\psi(N)=\frac1N\sum_{n=1}^N u_n.
\]
Since \(u_n\downarrow0\), we have \(\psi(N)\to0\) and
\(\psi(N)\ge u_N\).
With \(h_N=N\psi(N)\), \(h_0=0\),  
\[
0<h_N-h_{N-1}=u_N\le u_1.
\]
Moreover,
\[
\frac{\rho_N}{\psi(N)}\le\sqrt{\rho_N}\to0,
\qquad
\frac{\sqrt{(\log N)/N}}{\psi(N)}
\le N^{-1/4}\sqrt{\log N}\to0.
\]
Theorem~\ref{thm:prescribed-above-square-root} provides  a PPC sequence
with \(D_N^*\sim c_\psi\psi(N)\), and hence
\(D_N^*/\rho_N\to\infty\).
\end{proof}

In the proof of the main theorems, the low-discrepancy construction uses nested random point sets with
equal numbers of points in dyadic cells at each completed stage.
An ordering of the newly added points controls discrepancy for every
$N$, while conditional pair-correlation estimates give PPC almost
surely. A conditional-expectation argument, based on the method of
\citet{Raghavan1988}, selects one realization deterministically.
To obtain other discrepancy rates, we either insert sparse additional
points into a fixed sequence satisfying
Theorem~\ref{thm:low-discrepancy-background}, or sample independently
from suitable nonidentical distributions.

\medskip

Sections~\ref{sec:low-discrepancy-background}--\ref{sec:prescribed-bias}
prove Theorems~\ref{thm:low-discrepancy-background}--\ref{thm:prescribed-above-square-root},
respectively. Section~\ref{sec:consequences-scope} discusses optimal
discrepancy constants.

\section{Construction of a low-discrepancy PPC sequence}
\label{sec:low-discrepancy-background}

We prove Theorem~\ref{thm:low-discrepancy-background} in two stages. First,
we construct a random sequence every realization of which has low discrepancy,
and establish summable bounds on the probability that its pair-correlation
estimates deviate from their expected values. Second, we convert this
randomized construction into a deterministic one by means of a
conditional-expectation rule, thereby producing the deterministic sequence
required for Theorem~\ref{thm:low-discrepancy-background}.

\subsection{Nested random point sets}  
We describe a recursive scheme that generates a sequence of increasingly refined random point sets on
$\mathbb{T}$. 

\medskip

For \(k\ge1\), a dyadic interval of order \(k\), or a level-\(k\) cell, takes the form
\[
  I_{k,r}=\left[\frac{r}{2^k},\frac{r+1}{2^k}\right) 
  \qquad (0\le r<2^k).
\]
To each such interval we assign \(k\) labels \((k,r,j)\), \(1\le j\le k\).
Ordering all these labels lexicographically defines a bijection
\[
  \tau:\{(k,r,j): k\ge1,\ 0\le r<2^k,\ 1\le j\le k\}
  \to  \mathbb N,
\]
where \(\tau(k,r,j)=n\) means that \((k,r,j)\) is the \(n\)-th label in this
order.

We work on the product probability space
\[
  (\Omega,\mathcal F,\PP)
  =
  \bigl([0,1)^{\mathbb N},
  \mathcal B([0,1))^{\otimes\mathbb N},
  \lambda^{\otimes\mathbb N}\bigr).
\]
Define a sequence of random variables \((U_n)\) by
\[
  U_n(\omega)=\omega_n,
  \qquad \omega=(\omega_j)_{j\ge1}\in\Omega.
\]
Then the \(U_n\) are independent and uniform on \([0,1)\).
With the \(n\)-th label \((k,r,j)\), where \(n=\tau(k,r,j)\), we associate the
candidate
\[
  Z_n=\frac{r+U_n}{2^k},
\]
which is uniform on its assigned cell \(I_{k,r}\). Since  \(U_n\)
are   independent, the candidates \(Z_n\) are   independent as
well.

For each   \(\omega\in\Omega\), we now construct \(P_k(\omega)\) recursively. In our construction, all point counts are taken with multiplicity. We start
with
\[
  P_1(\omega)=\{Z_1(\omega),Z_2(\omega)\}.
\]
Assume that \(P_k(\omega)\) has been defined and contains exactly \(k\)
points in each \(I_{k,r}\). For \(0\le r<2^{k+1}\), put
\[
  p_r(\omega)=\sharp\bigl(P_k(\omega)\cap I_{k+1,r}\bigr),
  \qquad
  q_r(\omega)=k+1-p_r(\omega).
\]
Since \(p_r(\omega)\in\{0,\dots,k\}\), we have
\(q_r(\omega)\in\{1,\dots,k+1\}\). In each \((k+1)\)-cell
\(I_{k+1,r}\), we retain the first \(q_r(\omega)\) of the \(k+1\)
candidates, namely
\[
  Z_{\tau(k+1,r,j)}(\omega),
  \qquad 1\le j\le q_r(\omega).
\]
Adjoining these retained points to \(P_k(\omega)\), we obtain
\(P_{k+1}(\omega)\), which contains exactly \(k+1\) points in each
\(I_{k+1,r}\). Thus
\[
  N_k:=\sharp P_k(\omega)=k2^k,
  \qquad
  Q_k:=N_{k+1}-N_k
  =\sum_{r=0}^{2^{k+1}-1}q_r(\omega)=(k+2)2^k.
\] 

 Whether a candidate is retained is determined by the preceding point
counts. 
Each point count is a finite sum of indicators of the event that a candidate
is retained and belongs to a given set; induction therefore shows that all
retention events and point counts are measurable. Following the usual convention, we henceforth suppress the dependence on
\(\omega\). 

\begin{lemma}
\label{lem:lowdisc-completed-independence}
In terms of point counts, the random point set \(P_k\) is distributionally
equivalent to the model in which \(k\) independent uniform points are sampled
in each level-\(k\) cell, independently across cells. More precisely, for every
finite family of Borel subsets of \([0,1)\), the counts of \(P_k\) in these
sets have the same joint distribution as the corresponding counts obtained by
sampling \(k\) independent uniform points in each \(I_{k,r}\), independently across cells.
\end{lemma}

\begin{proof}
For this proof, we list the points in each cell as follows. For \(k=1\), each
level-\(1\) list consists of its single point. At each transition
\(k\to k+1\), the old points in each \((k+1)\)-cell inherit their relative
order from the list in its parent node, and the retained candidates are
then appended in increasing order of their indices.

We prove by induction on
\(k\) that all entries in the level-\(k\) lists are mutually independent and
uniform on their respective cells. 

The base case \(k=1\) follows from the
independence of the initial samples.

Assume the assertion at level \(k\). To pass to level \(k+1\), we condition  on the half containing every old point in
every level-\(k\) cell.   More precisely,
for each cell \(I_{k,j}\), let \(X_{j,1},\dots,X_{j,k}\) denote the old points
in this cell. Splitting \(I_{k,j}\) into its two offspring
\(I_{k+1,2j}\) and \(I_{k+1,2j+1}\), we condition on the joint event that each
\(X_{j,i}\) lies in a prescribed half, i.e., on the entire vector of
membership indicators
\(
  \bigl(
    \mathbf 1_{\{X_{j,i}\in I_{k+1,2j}\}}
  \bigr)_{0\le j<2^k,\ 1\le i\le k}.
\)
 By the induction hypothesis, the old points remain conditionally independent and uniform on their specified halves.

Under this conditioning, the counts \(p_r\) and 
\(q_r\) are fixed. The list in each \(I_{k+1,r}\) contains \(p_r\) old points followed by
the first \(q_r\) new candidates. The new candidates remain independent
of one another and of the old points.
All entries in the completed lists are therefore conditionally
independent and uniform on their respective cells.

For arbitrary Borel  sets \(H_{r,i}\subseteq I_{k+1,r}\),
\(0\le r<2^{k+1}\), \(1\le i\le k+1\), the conditional probability
that the \(i\)-th point in \(I_{k+1,r}\) belongs to \(H_{r,i}\)
for every \(r,i\) is
\(
  \prod_{r,i}
  {\lambda(H_{r,i})}/{\lambda(I_{k+1,r})}.
\)
This expression does not depend on the specified halves.
By the law of total probability, the same formula holds
without conditioning. This completes the induction and hence establishes the asserted equality of
joint count distributions.
\end{proof}

\subsection{Ordering the random points}

The   construction described above produces the point sets \(P_k\) but
does not specify the order in which their points are arranged. To obtain a
low-discrepancy sequence \(Y_1, Y_2, Y_3, \dots\), we order these points
recursively.

We set
\(Y_1=Z_1\) and \(Y_2=Z_2\) to be the two points in \(P_1\).
Suppose that the points of \(P_k\) have been ordered as
\(Y_1,\ldots,Y_{N_k}\).
We retain this order and specify the order in which the
\(Q_k\) new points are appended.

Consider the dyadic tree with root \(I_{0,0}=[0,1)\), 
and
leaves \(I_{k+1,0},\ldots,I_{k+1,2^{k+1}-1}\) at depth \(k+1\).
For each node \(J\), we define its quota by
\(
  q(J)=\sum_{I_{k+1,r}\subseteq J}q_r
\), and thus \(q(J)\) is the total number of new points added to \(J\) during the
current transition. In particular,
\(q(I_{k+1,r})=q_r\) and \(q\bigl([0,1)\bigr)=Q_k\).

At an internal node \(J\), we write \(J_{\mathrm L}\) and \(J_{\mathrm R}\)
for its left and right offspring, respectively.
We associate with \(J\) the (upper)  mechanical word  \(\chi_J =(\chi_J(1),\dots,\chi_J(q(J)))
  \in\{0,1\}^{q(J)}
\)
of slope \(q(J_{\mathrm R})/q(J)\), characterized by
\[
  \chi_J(n)=
  \left\lceil\frac{nq(J_{\mathrm R})}{q(J)}\right\rceil
  -
  \left\lceil\frac{(n-1)q(J_{\mathrm R})}{q(J)}\right\rceil 
  \qquad (1\le n\le q(J)).
\]  

In the following,  we issue \(Q_k\) requests in
total; each request starts at the root, descends along the tree, and
eventually lands in some leaf. Exactly \(q_r\) requests land in \(I_{k+1,r}\). The
order in which the requests reach the leaves will determine the order in which the new
points are placed.

 Every request starts at the root. The \(v\)-th request reads the \(v\)-th
        symbol \(\chi_{I_{0,0}}(v)\) of the   mechanical word  \(\chi_{I_{0,0}}\): if it is \(0\), the request
        descends into the left subtree; if it is \(1\), into the right
        subtree. 
        Suppose a request reaches a node
        \(J'\). That node records that it has now been visited \(v'\) times,
        and the \(v'\)-th symbol of \(\chi_{J'}\) determines whether the
        request descends into the left or the right offspring of \(J'\).
         The descent is repeated at every
        level until the request reaches some leaf \(I_{k+1,r}\), whence the request is
        assigned the label \(r\).

Reading the labels in the order
        \(1,2,\dots,Q_k\) of the requests yields
        \[
          W_k=(\ell_1,\ell_2,\dots,\ell_{Q_k}) 
          \qquad (\ell_n\in\{0,\dots,2^{k+1}-1\}).
        \]

The quota of the leaf \(I_{k+1,r}\) is \(q_r\), meaning that \(q_r\) requests should
land in \(I_{k+1,r}\). This is guaranteed by the mechanical word at each node: in
\(\chi_J\), the number of zeros equals the left quota \(q(J_L)\) and the
number of ones equals the right quota \(q(J_R)\), so among the requests
entering \(J\), exactly \(q(J_L)\) descend left and \(q(J_R)\) descend right.
Applying this level by level, the number of requests reaching each leaf equals
the quota of that leaf. Hence
\[
  \sharp\{n: \ell_n=r\}=q_r,
  \qquad
  \sum_{r=0}^{2^{k+1}-1} q_r=Q_k .
\]

To determine the order in which the points are finally appended, we need to
impose a definite order on the \(q_r\) points that the leaf \(I_{k+1,r}\)
contributes to \(P_{k+1}\).  Since these \(q_r\) points arise by retaining the first \(q_r\) elements of a larger sample, they already carry a natural order, namely the index order of
the sample. The internal ordering of a leaf \(I_{k+1,r}\) is taken to be precisely this index order.

We are ready to append the \(Q_k\) new points to the existing sequence
\(Y_1,\dots,Y_{N_k}\). Scan \(W_k\) in the order \(n=1,2,\dots,Q_k\). Suppose
\(\ell_n=r\), and suppose this is the \(w\)-th occurrence of the label \(r\) in
\(W_k\). Then the \(w\)-th point of the ordered list carried by
\(I_{k+1,r}\cap(P_{k+1}\setminus P_k)\) is appended to the sequence.

We have thus defined the order of the points in \(P_{k+1}\), which yields the
sequence
\[
  Y_1, Y_2, \dots, Y_{N_{k+1}}.
\]
We remark that each \(Y_n\) depends on only finitely many of the candidates
\(Z_m\), and is therefore a measurable function of those candidates.
The random variables \(Y_n\) need not be independent.

\begin{lemma}
\label{lem:lowdisc-binary-balance}
For every internal node \(J\) and an integer \(0\le v\le q(J)\),
\[
  \sharp\{1\le n\le v:\chi_J(n)=0\}
  =
  \left\lfloor\frac{vq(J_{\mathrm L})}{q(J)}\right\rfloor.
\]
Consequently,
\[
  \left|
  \sharp\{1\le n\le v:\chi_J(n)=0\}
  -v\frac{q(J_{\mathrm L})}{q(J)}
  \right|<1.
\]
The corresponding errors for the two offspring are negatives of each other.
\end{lemma}
\begin{proof}This is the balance property of mechanical words; we refer the reader to
\cite{Lothaire2002} for further details. For the reader's convenience, we
include a short proof.

We readily check that \(1-\chi_J(n)=
  \left\lfloor {nq(J_{\mathrm L})}/{q(J)}\right\rfloor
  -
  \left\lfloor {(n-1)q(J_{\mathrm L})}/{q(J)}\right\rfloor\), and thus  
\[
   \sharp\{1\le n\le v:\chi_J(n)=0\}=\sum_{n=1}^{v}
  \left(
  \left\lfloor\frac{nq(J_{\mathrm L})}{q(J)}\right\rfloor
  -
  \left\lfloor\frac{(n-1)q(J_{\mathrm L})}{q(J)}\right\rfloor
  \right)
  =
  \left\lfloor\frac{vq(J_{\mathrm L})}{q(J)}\right\rfloor.
\]
The remaining assertions follow from the floor bound and the fact
that the counts for the two offspring sum to \(v\).
\end{proof}

When no confusion can arise, we identify a label \(r\in\{0,\dots,2^{k+1}-1\}\)
with the leaf \(I_{k+1,r}\) it represents, and a set \(S\) of labels with the
union \(\bigcup_{r\in S}I_{k+1,r}\) of the corresponding leaves. Thus
\(q(S)=\sum_{r\in S}q_r\); we also write \(p(S)=\sum_{r\in S}p_r\).
For \(0\le m\le Q_k\) and \(0\le r<2^{k+1}\), we put
\[
  x_r(m)=\sharp\{1\le n\le m:\ell_n=r\},
  \qquad
  E_r(m)=x_r(m)-m\frac{q_r}{Q_k}.
\]


   \begin{lemma}
\label{lem:lowdisc-spatial-prefix}
For integers \(m,b\) with \(0\le m\le Q_k\) and \(0< b\le 2^{k+1}\), we have
\[
  \left|\sum_{r=0}^{b-1}E_r(m)\right|\le k+1.
\]
Consequently, for  \(0\le b_1< b_2\le 2^{k+1}\),
\(
  \left|\sum_{r=b_1}^{b_2-1}E_r(m)\right|\le2(k+1).
\)
\end{lemma}
 \begin{proof}
Throughout, requests are those issued during the current
transition \(k\to k+1\).

For a node \(J\) of the dyadic tree introduced earlier, let \(h(J)\) be the
height of \(J\), that is, the number of edges on a path from \(J\) to a leaf,
so that \(h(J)=0\) at a leaf and \(h(I_{0,0})=k+1\). 

We consider only the requests reaching \(J\). Let \(S\) be an initial
segment of the leaf labels below \(J\); thus \(S\) is either empty, or all
labels below \(J\), or a set of the form \(\{0,\dots,b-1\}\) intersected with
the labels below \(J\).
For \(0\le v\le q(J)\),  let
\(F_J(v,S)\) be the number of the first \(v\) such requests that end in a leaf
labeled by \(S\), minus the ideal count \(v\,q(S)/q(J)\):
\[
  F_J(v,S)
  =
  \sharp\{1\le n\le v:\text{the }n\text{-th request reaching }J
  \text{ ends in a leaf of }S\}
  -
  v\frac{q(S)}{q(J)}.
\] 
We claim that, for every node \(J\), every
\(0\le v\le q(J)\), and every initial segment \(S\) of the labels below
\(J\),
\[
  |F_J(v,S)|\le h(J).
\]

We argue by induction on \(h(J)\). When \(h(J)=0\), i.e., \(J\) is a leaf, the set \(S\) is either empty or
consists of the single label of \(J\). In either case, \(F_J(v,S)=0\).

Let \(J\) be a  node with offspring \(J_{\mathrm L}\) and
\(J_{\mathrm R}\). Let \(v_0\) and \(v_1\) be the numbers of the first \(v\)
requests reaching \(J\) that descend into \(J_{\mathrm L}\) and
\(J_{\mathrm R}\), respectively, so \(v_0+v_1=v\).   Set
\(
  \delta=v_0-v {q(J_{\mathrm L})}/{q(J)}.
\)
By the balance property established in Lemma \ref{lem:lowdisc-binary-balance}, we have
\(|\delta|<1\). 
We distinguish two cases, according to the position of \(S\).

\smallskip

\emph{Case 1: all leaves labelled by \(S\) lie below \(J_{\mathrm L}\).} Then a
request reaching \(J\) ends in \(S\) if and only if it descends into
\(J_{\mathrm L}\) and ends in \(S\) there. So
\[
  F_J(v,S)
  =
  F_{J_{\mathrm L}}(v_0,S)
  +
  v_0\frac{q(S)}{q(J_{\mathrm L})}
  -
  v\frac{q(S)}{q(J)}
  =
  F_{J_{\mathrm L}}(v_0,S)
  +
  \delta\frac{q(S)}{q(J_{\mathrm L})}.
\]

\smallskip

\emph{Case 2: \(S\) contains all leaves below \(J_{\mathrm L}\).} Write
\(S=J_{\mathrm L}\cup S_1\), where \(S_1\) is an initial segment of the labels
below \(J_{\mathrm R}\). A request reaching \(J\) ends in
\(S\) if and only if it either descends into \(J_{\mathrm L}\), or descends
into \(J_{\mathrm R}\) and ends in \(S_1\). Hence
\[
  F_J(v,S)
  =
  v_0
  +
  F_{J_{\mathrm R}}(v_1,S_1)
  +
  v_1\frac{q(S_1)}{q(J_{\mathrm R})}
  -
  v\frac{q(J_{\mathrm L})+q({S_1})}{q(J)}.
\]
Using \(v_1=v-v_0\) and the definition of \(\delta\), this simplifies to
\[ 
  F_J(v,S)
  =
  F_{J_{\mathrm R}}(v_1,S_1)
  +
  \delta\left(1-\frac{q(S_1)}{q(J_{\mathrm R})}\right).
\]

In both cases, the induction hypothesis applied to the offspring gives
\[
  |F_J(v,S)|
  \le
  \bigl(h(J)-1\bigr)+|\delta|
  \le
  h(J),
\]
which proves the claim at \(J\).

Now taking 
\(S=\{0,1,\dots,b-1\}\), the labels below the root \(I_{0,0}\) are all leaf labels,  so
\[
  F_{I_{0,0}}(m,S)
  =
  \sum_{r=0}^{b-1}x_r(m)
  -
  m\frac{\sum_{r=0}^{b-1}q_r}{Q_k}
  =
  \sum_{r=0}^{b-1}E_r(m).
\]
The claim at the root therefore yields
\[
  \left|\sum_{r=0}^{b-1}E_r(m)\right|
  \le h(I_{0,0})=k+1,
\]
which is the first estimate. The second then follows from the triangle
inequality. \end{proof}
 
\subsection{Discrepancy of every prefix}

We calculate the discrepancy of the sequence. To this end, we fix a transition from level $k$ to level $k+1$. 
 Recall that \(N_k=k2^k\) and \(Q_k=(k+2)2^k\) are the numbers of old points
and new points in the transition, respectively.

\begin{lemma}
\label{lem:lowdisc-completed-spatial-count}
For an integer $b$ with  $0\le b\le 2^{k+1}$, writing \(S_b=\{0,\ldots,b-1\}\), we have
\begin{equation}
\label{eq:lowdisc-completed-spatial-count}
\left|p (S_b) -N_k\frac{b}{2^{k+1}}\right|\le\frac{k}{2},
\end{equation}
and
\begin{equation}
\label{eq:lowdisc-quota-spatial-count}
q ({S_b})-Q_k\frac{b}{2^{k+1}} =-\left(p({S_b})-N_k\frac{b}{2^{k+1}} \right).
\end{equation}
\end{lemma}

\begin{proof}
Each level-$k$ cell contains $k$ old points, so
\[p (S_b) -N_k\frac{b}{2^{k+1}}=
\begin{cases}
0,&b=2j,\\
p_{2j}-k/2,&b=2j+1.
\end{cases}
\]
Since \(0\le p_{2j}\le k\), this yields
\eqref{eq:lowdisc-completed-spatial-count};
\eqref{eq:lowdisc-quota-spatial-count} is equivalent to the identity
\(p(S_b)+q(S_b)=(k+1)b\). 
\end{proof}

\begin{proposition}
\label{prop:lowdisc-prefix-discrepancy}
The sequence \(Y_1,Y_2,\ldots\) constructed above satisfies
\begin{equation}
\notag
ND_N^*\le\frac{5}{2\log2}\log N
\qquad(N\ge2).
\end{equation}
\end{proposition}

\begin{proof}
Write \(N=N_k+m\) with \(k\ge1\) and \(0\le m\le Q_k\). Throughout this
section, every integer \(N\) is assumed to take this form, so that \(N\) and
the pair \((k,m)\) determine each other.

For an integer $b$ with $0\le b\le 2^{k+1}$,  
\eqref{eq:lowdisc-quota-spatial-count} yields
\begin{align*}
\sharp\{1\le n\le N:Y_n<\frac{b}{2^{k+1}}\}-N\frac{b}{2^{k+1}}&=p ({S_b})+\sum_{r=0}^{b-1}x_r(m)-N\frac{b}{2^{k+1}}\\
&=\left(1-\frac{m}{Q_k}\right)\left(p ({S_b})-N_k\frac{b}{2^{k+1}}\right)
 +\sum_{r=0}^{b-1}E_r(m).    
\end{align*}
By Lemma~\ref{lem:lowdisc-spatial-prefix} and
\eqref{eq:lowdisc-completed-spatial-count},
\begin{equation}
\label{eq:lowdisc-boundary-error}
\left|\sharp\{1\le n\le N:Y_n<\frac{b}{2^{k+1}}\}-N\frac{b}{2^{k+1}}\right|
\le\left(1-\frac{m}{Q_k}\right)\frac{k}{2}+k+1
\le\frac{3k}{2}+1.
\end{equation}

For $0\le t<1$, choose an integer $b$ so that $ {b}/{2^{k+1}}\le t<{(b+1)}/{2^{k+1}}$.
Since each level-$(k+1)$ cell contains at most $k+1$ of the first $N$ points, and
$N\le(k+1)2^{k+1}$, both
\(
\sharp\{1\le n\le N: {b}/{2^{k+1}}\le Y_n<t\}\) and \(N(t-{b}/{2^{k+1}})\) 
lie in $[0,k+1]$, so their difference has absolute value at most $k+1$.
Together with \eqref{eq:lowdisc-boundary-error} and the exact count at
$t=1$, this yields
\(
ND_N^*\le {5k}/{2}+2.
\)

For $N\ge8$, we have $k\ge2$ and $N\ge k2^k$, and thus
\[
\frac{5k}{2}+2
\le\frac{5k}{2}+\frac{5\log k}{2\log2}
=\frac{5}{2\log2}\log(k2^k)
\le\frac{5}{2\log2}\log N.
\]

For $2\le N\le7$, the bound
$ND_N^*\le N\le5\log N/(2\log2)$ completes the proof.
\end{proof}

\subsection{Conditional pair-correlation estimates}

We  estimate the conditional expectation of \(NR_N(s)\).

\begin{lemma}
\label{lem:lowdisc-scheduling-energy}
For \(k\ge1\) and  \(0\le m\le Q_k\), we have
\begin{equation}
\notag
  \sum_{r=0}^{2^{k+1}-1}E_r(m)^2\le 5\cdot2^k.
\end{equation}
\end{lemma}
\begin{proof}
 For \(0\le j<2^k\), let \(v_j(m)\) be the number
of visits to \(I_{k,j}\) among the first \(m\) requests.
Since \(q_{2j}+q_{2j+1}=k+2\) for every \(j\), requests
alternate between the two offspring at each node of level
less than \(k\), giving
\[
  v_j(m)\in
  \left\{
    \left\lfloor\frac{m}{2^k}\right\rfloor,
    \left\lceil\frac{m}{2^k}\right\rceil
  \right\},
  \qquad
  \left|v_j(m)-\frac{m}{2^k}\right|
  <1.
\]
 By Lemma~\ref{lem:lowdisc-binary-balance},  for \(r\in\{2j,2j+1\}\), we have that \(|x_r(m)-v_j(m){q_r}/({k+2})|<1\), and thus
\begin{align*}
  |E_r(m)|
  \le
  \left|x_r(m)-\frac{q_r}{k+2}v_j(m)\right|
  +\frac{q_r}{k+2}\left|v_j(m)-\frac{m}{2^k}\right|
  <1+\frac{q_r}{k+2}.
\end{align*}
Hence
\begin{align*}
  E_{2j}(m)^2+E_{2j+1}(m)^2
  <\left(1+\frac{q_{2j}}{k+2}\right)^2
    +\left(1+\frac{q_{2j+1}}{k+2}\right)^2\le 5.
\end{align*}
Summing over the \(2^k\) pairs completes the proof.
\end{proof}

By Lemma~\ref{lem:lowdisc-completed-independence},
the variables \(p_{2j}\), \(0\le j<2^k\), are independent, with
\[
  p_{2j}\sim\operatorname{Binomial}(k,1/2),
  \qquad p_{2j+1}=k-p_{2j}.
\]
For   \(0\le m\le Q_k\), recalling \(N=N_k+m\), we define, for \(0\le r<2^{k+1}\),
\[
  C_r(m)=p_r+x_r(m),
  \qquad
  e_r(m)=C_r(m)-\frac{N}{2^{k+1}},
\]
where \(C_r(m)\) is the number of points in \(I_{k+1,r}\) among the first
\(N\) terms and \(e_r(m)\) is its deviation from \(N/2^{k+1}\).
Define the event
\begin{equation}
\notag
  \mathcal A_k=
  \left\{
    \sum_{r=0}^{2^{k+1}-1}\left(p_r-\frac{k}{2}\right)^2\le k2^{k+1}
  \right\}.
\end{equation}

\begin{lemma}
\label{lem:lowdisc-current-energy}
For \(k\ge1\), we have
\begin{equation}
\label{eq:lowdisc-occupation-failure}
  \PP(\mathcal A_k^c)\le\frac{1}{9\cdot2^{k-1}}.
\end{equation}
On \(\mathcal A_k\),  we have, uniformly in \(0\le m\le Q_k\),
\begin{equation}
\notag
  \sum_{r=0}^{2^{k+1}-1}e_r(m)^2\le 	k2^{k+4}.
\end{equation}
\end{lemma}

\begin{proof}
Set
\[
  X=\sum_{j=0}^{2^k-1}\left(p_{2j}-\frac{k}{2}\right)^2.
\]

Since $p_{2j}\sim\operatorname{Binomial}(k,1/2)$, a direct calculation yields
\[
    \E\!\left[\left(p_{2j}-\frac{k}{2}\right)^{2}\right]=\frac{k}{4},
    \qquad
    \E\!\left[\left(p_{2j}-\frac{k}{2}\right)^{4}\right]
    =\frac{3k^{2}}{16}-\frac{k}{8},
\]
and hence
\[
    \operatorname{Var}\!\left(\left(p_{2j}-\frac{k}{2}\right)^{2}\right)
    =\frac{k(k-1)}{8}\le\frac{k^{2}}{8}.
\]
The  independence property then gives $\E X=k2^{k-2}$ and $\operatorname{Var}(X)\le k^{2}2^{k-3}$.
Since \(p_{2j+1}=k-p_{2j}\), we see that 
\(
  \sum_{r=0}^{2^{k+1}-1}\left(p_r-{k}/{2}\right)^2=2X
\),  and thus
  \(X-\E X>3k2^{k-2}\) on \(\mathcal A_k^c\).   Chebyshev's
inequality then yields
\[
  \PP(\mathcal A_k^c)
  \le
  \frac{k^2 2^{k-3}}{(3k2^{k-2})^2}
  =\frac{1}{9\cdot2^{k-1}}.
\]

As in Lemma \ref{lem:lowdisc-completed-spatial-count}, we have \(
  q_r-{Q_k}/{2^{k+1}}=-\left(p_r-{k}/{2}\right).
\)
Therefore, by the  definitions of \(e_r(m)\) and \(E_r(m)\), we obtain 
\begin{equation}
\notag
  e_r(m)=
  \left(1-\frac{m}{Q_k}\right)
  \left(p_r-\frac{k}{2}\right)+E_r(m).
\end{equation}
For  \(0\le m\le Q_k\),  Lemma~\ref{lem:lowdisc-scheduling-energy}
and the elementary inequality
\( (a+b)^2\le 2(a^2+b^2)\)
  imply,  on  \(\mathcal A_k\), 
\[
  \sum_{r=0}^{2^{k+1}-1} e_r(m)^2
  \le 2(k2^{k+1}+5\cdot 2^k)
 < k2^{k+4}.
\]
\vskip -0.88cm
\end{proof}

\medskip

Let
\(
\mathcal G_k=\sigma(p_0,\ldots,p_{2^{k+1}-1})\subseteq\mathcal F
\)
be the $\sigma$-algebra generated by the old-point counts in
$I_{k+1,0},\ldots,I_{k+1,2^{k+1}-1}$.
These counts determine $q_r=k+1-p_r$ and the order of the new points.
For each integer \(0\le m\le Q_k\), conditioning on \(\mathcal G_k\) therefore
fixes, for each of the first \(m\) new points, both the cell containing it and
its candidate index, and hence the counts \(x_r(m)\) and \(C_r(m)\).

\begin{lemma}
\label{lem:lowdisc-conditional-independence}
Fix  \(0\le m\le Q_k\), and   \(N=N_k+m\).
Conditional on \(\mathcal G_k\), the counts of \(Y_1,\ldots,Y_N\)
in any finite family of Borel subsets of \([0,1)\) have the
same joint distribution as the corresponding counts obtained
by sampling \(C_r(m)\) independent uniform points in each
\(I_{k+1,r}\), independently across cells.
\end{lemma}

\begin{proof}
By Lemma~\ref{lem:lowdisc-completed-independence} and the
conditioning argument in its proof, conditional on
\(\mathcal G_k\), the joint count distribution of the old
points is that obtained from \(p_r\) independent uniform
points in each \(I_{k+1,r}\), independently across cells.

The new candidates are independent of all old points.
Given \(\mathcal G_k\), the first \(m\) new points include
the first \(x_r(m)\) candidates in \(I_{k+1,r}\), with
\(x_r(m)\) fixed.
Adjoining these candidates proves the assertion, since
\(p_r+x_r(m)=C_r(m)\).
\end{proof}

Fix $s>0$. Since $N=N_k+m\ge k2^k$, 
\begin{equation}
\notag
a:=\frac{s/N}{1/2^{k+1}}
\le\frac{2s}{k}
\to 0
\qquad (N\to\infty).
\end{equation}
Hence, for all sufficiently large \(N\), every pair counted by \(R_N(s)\)
lies either in the same cell or in two adjacent cells on \(\mathbb T\). In the
latter case we say that the pair crosses a boundary, by which we mean the
common endpoint shared by two adjacent cells on \(\mathbb T\); there are
\(2^{k+1}\) such boundaries.

\begin{lemma}
Whenever $a<1$, we have
\begin{align}
\E\left[NR_N(s)\mid\mathcal G_k\right]
&=
(2a-a^2)
\sum_{r=0}^{2^{k+1}-1}C_r(m)\bigl(C_r(m)-1\bigr)
\notag\\
&\quad+
a^2
\sum_{r=0}^{2^{k+1}-1}C_r(m)C_{r+1}(m),
\label{eq:lowdisc-conditional-expectation}
\end{align}
where the cell indices are interpreted modulo $2^{k+1}$.
\end{lemma}

\begin{proof}
Since $NR_N(s)$ is unchanged by permuting the first $N$ points,
Lemma~\ref{lem:lowdisc-conditional-independence} allows us to
compute its conditional expectation using independent uniform
samples, with $C_r(m)$ points in each $I_{k+1,r}$.

For two independent uniform points in the same cell, rescaling the cell to
\([0,1)\) gives the probability
\(
  2a-a^2.
\)
There are $C_r(m)(C_r(m)-1)$ ordered pairs of terms in $I_{k+1,r}$
with different indices, giving the first term on the right-hand side of
\eqref{eq:lowdisc-conditional-expectation}.

For an ordered pair that crosses a boundary, we interchange the two points if
necessary, so that the first lies in the left cell and the second in the right
cell, and then rescale the two adjacent cells to \([0,1)\) and \([1,2)\). We thus see that  the probability that  the   pair satisfies the distance condition is
\(
 {a^2}/{2}.
\)
There are $2C_r(m)C_{r+1}(m)$ ordered pairs across a boundary,
so their conditional expected contribution is
\[
2C_r(m)C_{r+1}(m)\frac{a^2}{2}
=a^2C_r(m)C_{r+1}(m).
\]
Summing over the $2^{k+1}$  boundaries gives the second term.
\end{proof}

\subsection{Almost-sure Poissonian pair correlation}

We prove that $(Y_n(\omega))_{n\ge1}$ has PPC for
$\PP$-almost every $\omega\in\Omega$.

Lemma~\ref{lem:lowdisc-conditional-independence} identifies the conditional
distribution of the first \(N\) points given \(\mathcal G_k\) with the
distribution obtained by sampling \(C_r(m)\) independent uniform points in each
\(I_{k+1,r}\), independently across cells.
Since $s/N\le2^{-(k+2)}$, moving one point within its cell can affect
only pairs containing that point and another point in the same
cell or an adjacent cell.
Each cell contains at most $k+1$ points, so at most $6(k+1)$
ordered-pair indicators can change.
Thus changing one point changes $R_N(s)$ by at most $6(k+1)/N$.
Applying McDiarmid's inequality \cite{McDiarmid1989} to these
independent samples and using the equality of conditional
distributions gives, for every $\varepsilon>0$,
\begin{equation}
\label{eq:lowdisc-pair-concentration}
\PP\bigl(
\left|R_N(s)-\E[R_N(s)\mid\mathcal G_k]\bigr|>\varepsilon
\,\middle|\,\mathcal G_k
\right)
\le
2\exp\left(-\frac{\varepsilon^2N}{18(k+1)^2}\right).
\end{equation}
Taking expectations gives the same bound for the unconditional probability. Therefore
\begin{align*}
&\sum_{k\ge4s}\sum_{N=N_k}^{N_{k+1}}
\PP\left(
\left|R_N(s)-\E[R_N(s)\mid\mathcal G_k]\right|>\varepsilon
\right)\\
 \le ~&
2\sum_{k\ge4s}(Q_k+1)
\exp\left(-\frac{\varepsilon^2N_k}{18(k+1)^2}\right)
<\infty.
\end{align*}
By the Borel--Cantelli lemma, we deduce that
\begin{equation}
    \label{limit1}
     R_N(s)-\E[R_N(s)\mid\mathcal G_k]\to 0
  \quad\text{almost surely as }N\to\infty,
\end{equation}
where \(k=k(N)\) is the unique integer with \(N_k\le N<N_{k+1}\).
Moreover, \eqref{eq:lowdisc-occupation-failure} gives
\(\sum_k\PP(\mathcal A_k^c)<\infty\). Applying the Borel--Cantelli lemma
once more, we obtain that \(\mathcal A_k\) holds for all sufficiently large
\(k\), almost surely.  

Substituting $C_r(m)=N/2^{k+1}+e_r(m)$ into
\eqref{eq:lowdisc-conditional-expectation}, and
using $\sum_r e_r(m)=0$ and the Cauchy--Schwarz inequality
\(
\left|\sum_r e_r(m)e_{r+1}(m)\right|
\le\sum_r e_r(m)^2,
\)
we obtain that, on  \(\mathcal A_k\), 
\begin{equation}
\label{eq:lowdisc-conditional-mean-bound}
\Bigl|\E[R_N(s)\mid\mathcal G_k]-2s\Bigr|
\le 2a\left(1+\frac1N\sum_r e_r(m)^2\right) 
\le34a\le\frac{68s}{k}.
\end{equation}
The above estimate uses \(0<a\le 2s/k\le 1/2\) and
\(N^{-1}\sum_r e_r(m)^2\le 16\). The second bound follows from
Lemma~\ref{lem:lowdisc-current-energy} together with \(N\ge N_k\).

Since \(\mathcal A_k\) holds almost surely for all sufficiently large \(k\),
and since \eqref{eq:lowdisc-conditional-mean-bound} holds on \(\mathcal A_k\),
we conclude that almost surely, 
\[
  \Bigl|\E[R_N(s)\mid\mathcal G_k]-2s\Bigr|
  \le\frac{68s}{k}\to 0.
\]  
Combining this with \eqref{limit1} yields   \(R_N(s)\to2s\) almost surely. Intersecting
the corresponding probability-one events over the countable set
\(s\in\mathbb Q_{>0}\), and using the monotonicity of \(R_N\) in \(s\) to
pass from rational to real \(s\), we obtain PPC almost surely.

\subsection{Deterministic selection}

We select the candidate points by a deterministic conditional-expectation
rule; see \citet{Raghavan1988} for   the method behind it. 

For $k\ge1$, put
\[
L_k=\lfloor\log_2 k\rfloor,
\qquad
\mathcal S_k=\left\{\frac pq:1\le p,q\le L_k,\ p,q\in\mathbb N\right\},
\qquad
\eta_k=k^{-\frac14},
\]
and define
\begin{equation}
\notag
\mathcal D_k
=
\bigcup_{\substack{ N_k\le N< N_{k+1}\\
s\in\mathcal S_k}}
\left\{|R_N(s)-2s|>\eta_k\right\}.
\end{equation}

Take \(k\) sufficiently large, say \(k\ge K_0\), that
\(L_k\le k/4\) and \(68L_k/k\le\eta_k/2\).
Then, by \eqref{eq:lowdisc-conditional-mean-bound}, for   $N_k\le N < N_{k+1}$ and $s\in\mathcal S_k$, we have, on $\mathcal A_k$,
\[
\left|\E[R_N(s)\mid\mathcal G_k]-2s\right|
\le\frac{\eta_k}{2}.
\]
Since $\mathcal A_k\in\mathcal G_k$, applying
\eqref{eq:lowdisc-pair-concentration} with $\varepsilon=\eta_k/2$,
multiplying by $\1_{\mathcal A_k}$, and taking expectations, we obtain 
\[
\PP\left(\mathcal A_k\cap\{|R_N(s)-2s|>\eta_k\}\right)
\le2\exp\left(-\frac{\eta_k^2N}{72(k+1)^2}\right).
\]
Splitting according to $\mathcal A_k$ and using the union bound,
we obtain
\begin{align}
\PP(\mathcal D_k)
&\le\PP(\mathcal A_k^c)
+\sum_{N=N_k}^{N_{k+1}-1}\sum_{s\in\mathcal S_k}
\PP\left(\mathcal A_k\cap\{|R_N(s)-2s|>\eta_k\}\right)
\notag\\
&\le\frac4{9\cdot2^{k+1}}
+2L_k^2Q_k
\exp\left(-\frac{\eta_k^2N_k}{72(k+1)^2}\right).
\notag
\end{align}
Since
$\eta_k^2N_k/(k+1)^2\asymp2^k/k^{3/2}$, the bound is summable in $k$.
Let $K_1\ge K_0$ be the least positive integer such that
\(\sum_{k\ge K_1}\PP(\mathcal D_k)<1/2.
 \)

For the candidate $Z_\nu$ assigned to $I_{b_\nu,r_\nu}$, write
\begin{equation}
\notag
Z_\nu
=
2^{-b_\nu}
\left(r_\nu+\sum_{j\ge1}\xi_{\nu,j}2^{-j}\right).
\end{equation}
The first \(b_\nu\) binary digits are fixed by the assigned cell, and
\(\xi_{\nu,j}\) is the \(j\)-th subsequent binary digit. We refer to all the
\(\xi_{\nu,j}\) as the free binary digits of \(Z_\nu\).
Under $\PP$, all the free digits $(\xi_{\nu,j})_{v,j\ge1}$ are independent and take
the values $0$ and $1$ with equal probability.
Arrange  these digits $(\xi_{\nu,j})_{v,j\ge1}$ in one sequence $(\zeta_m)_{m\ge1}$,
here we order the pairs $(\nu,j)$ by increasing $\nu+j$ and then
increasing $\nu$:
\[
(\zeta_1,\zeta_2,\ldots)
=
(\xi_{1,1},\xi_{1,2},\xi_{2,1},
 \xi_{1,3},\xi_{2,2},\xi_{3,1},\ldots).
\]

For $\nu,j\ge1$, define
\begin{equation}
\notag
\mathcal C_{\nu,j}
=
\{\xi_{\nu,j}=\xi_{\nu,j+1}=\cdots=\xi_{\nu,\,2j+\nu+1}\}.
\end{equation}
From the independence we obtain that 
\(\PP(\mathcal C_{\nu,j})=2^{-1-\nu-j}\) and \(\sum_{\nu,j\ge1}\PP(\mathcal C_{\nu,j})=1/2.\)
Define
\begin{equation}
\notag
\Phi
=
\sum_{k\ge K_1}\1_{\mathcal D_k}
+\sum_{\nu,j\ge1}\1_{\mathcal C_{\nu,j}}.
\end{equation}
The preceding estimates yield \(\E\Phi<1.\)

For a finite binary word $\sigma$ of length $n$, let
\[
[\sigma]=\{(\zeta_1,\ldots,\zeta_n)=\sigma\}, \quad \text{ with }~  [\varnothing]=\Omega.
\]
Then 
\(V(\sigma):=\E(\Phi\mid[\sigma])
=2^n\E[\Phi\,\1_{[\sigma]}]\), since \(\PP(\sigma)=2^{-n}\).
These conditional expectations are finite and satisfy
 \(V(\sigma)=(V(\sigma0)+V(\sigma1))/2.\)

Starting with $\sigma_0=\varnothing$, we define recursively
\begin{equation}
\notag
\sigma_n=
\begin{cases}
\sigma_{n-1}0,&V(\sigma_{n-1}0)\le V(\sigma_{n-1}1),\\
\sigma_{n-1}1,&V(\sigma_{n-1}1)<V(\sigma_{n-1}0).
\end{cases}
\end{equation}
We then obtain a nested sequence of words \((\sigma_n)\), and
\begin{equation}
\notag
  V(\sigma_n)\le V(\sigma_{n-1})\le V(\varnothing)<1
  \qquad(n\ge1).
\end{equation}
These words \((\sigma_n)\) determine a unique infinite digit sequence, and
the corresponding selected values of \(U_\nu\) thereby determine a point
\(\omega^*\in\Omega\).

None of the events \(\mathcal C_{\nu,j}\) can occur for the selected digit
sequence; otherwise, for \(n\) sufficiently large, we would have
\(\Phi\ge1\) on \([\sigma_n]\), and hence \(V(\sigma_n)\ge1\). Therefore no
candidate's free digits are eventually constant, and the corresponding
candidates are not dyadic rationals and lie in the interiors of their
assigned cells.

\medskip

Set \(y_n=Y_n(\omega^*)\). We show that \((y_n)\) has PPC by bounding
\(R_N^y(s)=R_N(s)(\omega^*)\). 

\smallskip

\def\sn{s^{{\scriptscriptstyle -}}} 
\def\sp{s^{{\scriptscriptstyle +}}} 

  Fix \(s>0\) and choose rational numbers \(\sn ,\sp \) with
\(
0<\sn <s<\sp .
\) Choose \(K_2\ge K_1\) such that
\(
\sn ,\sp \in\mathcal S_k
\) for \(k\ge K_1\). Fix \(k\ge K_1\) and  \(N\) satisfying
\(
N_k\le N<N_{k+1}.
\) Let \(\omega\in[\sigma_n]\), where \(\sigma_n\) is a selected prefix
whose length will be chosen below. We compare
\[
R_N(\sn )(\omega),\qquad R_N(s)(\omega^*),\qquad R_N(\sp )(\omega).
\]
These quantities depend only on \(
Y_1(\omega),  Y_2(\omega), \ldots, Y_N(\omega)\) and \(Y_1(\omega^*), Y_2(\omega^*), \ldots,Y_N(\omega^*).\)
The selection and ordering of the first \(N_{k+1}\) terms involve only
finitely many candidates, namely those up to level \(k+1\).

At \(\omega^*\), every candidate lies in the interior of its dyadic cell.
Let \(\delta>0\) be the minimum distance of their values from the
dyadic boundaries at level \(k+1\).
For a sufficiently long selected prefix \(\sigma_n\), we have 
\(
|Z_\nu(\omega)-Z_\nu(\omega^*)|<{\delta}/{2}\) for any \(\omega\in[\sigma_n]
\).  Hence \(Z_\nu(\omega)\) and
\(Z_\nu(\omega^*)\) lie in the same cell at each level up to \(k+1\).
Thus, for every \(\omega\in[\sigma_n]\), the first \(N\) terms use the
same candidate indices as at \(\omega^*\).

Extending the selected prefix \(\sigma_n\) further if necessary, we obtain
\[
  \|Y_i(\omega)-y_i\|_{\mathbb T}
  <
  \frac{\min\{s-\sn,\ \sp-s\}}{2N}
  \qquad(1\le i\le N,\ \omega\in[\sigma_n]).
\]
Consequently, every indexed pair at distance at most \(\sn/N\) for \(\omega\)
is at distance at most \(s/N\) for \(\omega^*\), and every indexed pair at
distance at most \(s/N\) for \(\omega^*\) is at distance at most \(\sp/N\)
for \(\omega\). Therefore
\[
  R_N(\sn)(\omega)\le R_N(s)(\omega^*)\le R_N(\sp)(\omega)
  \qquad(\omega\in[\sigma_n]).
\]

We claim that
\[
  2\sn-\eta_k\le R_N(s)(\omega^*)\le 2\sp+\eta_k.
\]
Indeed, if \(R_N(s)(\omega^*)>2\sp+\eta_k\), then
\(R_N(\sp)(\omega)>2\sp+\eta_k\) for every \(\omega\in[\sigma_n]\),
implying \([\sigma_n]\subseteq\mathcal D_k\) and hence
\(V(\sigma_m)\ge1\), a contradiction. The lower bound
\(R_N(s)(\omega^*)\ge 2\sn-\eta_k\) follows in the same way.

Letting \(N\to\infty\), and hence \(k\to\infty\), we obtain
\[
  2\sn
  \le \liminf_{N\to\infty}R_N(s)(\omega^*)
  \le \limsup_{N\to\infty}R_N(s)(\omega^*)
  \le 2\sp.
\]
Letting \(\sn,\sp\to s\), we conclude that
\(\lim_{N\to\infty}R_N(s)(\omega^*)=2s\), which proves that \((y_n)\) has PPC.
Proposition~\ref{prop:lowdisc-prefix-discrepancy} yields
\eqref{eq:main-low-discrepancy}. This completes the proof of
Theorem~\ref{thm:low-discrepancy-background}.

\section{PPC sequences with discrepancy of order
\texorpdfstring{$(\log N)/N$}{(log N)/N}}

Let \((y_n)_{n\ge1}\) be a sequence satisfying the conclusions of
Theorem~\ref{thm:low-discrepancy-background}. Given such a sequence, we prove
Theorem~\ref{thm:critical-order} by inserting zeros so as to obtain a matching
lower discrepancy bound for every sufficiently large \(N\), while preserving
PPC.

Set
\[
H_0=0,
\qquad
H_N=\lfloor\log N\rfloor,
\qquad
M_N=N-H_N
\quad(N\ge1).
\]
Then \(H_N-H_{N-1}\in\{0,1\}\) for every \(N\ge1\).
Retaining the original order of \((y_n)\), define
\begin{equation}
\notag
x_N=
\begin{cases}
0,&H_N-H_{N-1}=1,\\
y_{M_N},&H_N-H_{N-1}=0.
\end{cases}
\end{equation}
The first \(N\) terms thus consist of \(y_1,\ldots,y_{M_N}\)
and \(H_N\) copies of \(0\).

\begin{proposition}
The sequence \((x_N)\) has Poissonian pair correlation and satisfies
\begin{equation}
\label{eq:critical-two-sided}
\log N-1
\le ND_N^*
\le \left(1+\frac{5}{2\log2}\right)\log N
\qquad(N\ge2).
\end{equation}
In particular, \(D_N^*\asymp(\log N)/N\).
\end{proposition}
\begin{proof}
For $0<t\le1$, the construction gives
\begin{equation}
\notag
A_N^x(t)-Nt
=
\bigl(A_{M_N}^y(t)-M_Nt\bigr)+H_N(1-t).
\end{equation}
Since \(M_N\ge2\) for \(N\ge2\), we apply
\eqref{eq:main-low-discrepancy} and obtain
\begin{equation}
\notag
ND_N^*
\le H_N+\frac{5}{2\log2}\log M_N
\le \left(1+\frac{5}{2\log2}\right)\log N.
\end{equation}
For the lower bound, the first \(N\) terms contain at least \(H_N\) zeros.
Letting \(t\downarrow0\) in the definition of star discrepancy, we get
\[
ND_N^*\ge\sharp\{1\le n\le N:x_n=0\}
\ge H_N\ge\log N-1.
\]
This proves \eqref{eq:critical-two-sided}.

To verify PPC, fix \(s>0\). Splitting the ordered pairs according to
whether their entries are background points or inserted zeros, we obtain
\[
\begin{aligned}
R_N^x(s)
={}&\frac{M_N}{N}
R_{M_N}^y\left(s\frac{M_N}{N}\right)
+\frac{H_N(H_N-1)}{N}\\
&+\frac{2H_N}{N}
\sharp\left\{1\le m\le M_N:\normT{y_m}\le\frac{s}{N}\right\}.
\end{aligned}
\]
The first term tends to \(2s\) as \(N\to\infty\), by the PPC of
\((y_m)\), monotonicity in the window parameter, and \(M_N/N\to1\).
Applying \eqref{eq:main-low-discrepancy} to the two intervals adjacent
to \(0\), we obtain that, for all sufficiently large \(N\),
\[
\sharp\left\{1\le m\le M_N:\normT{y_m}\le\frac{s}{N}\right\}
\le \frac{2sM_N}{N}+\frac{5}{\log2}\log M_N
\ll_s\log N.
\]
Since \(H_N\le\log N\), the remaining two terms are
\(O_s((\log N)^2/N)\) and tend to \(0\).
Thus \(R_N^x(s)\to2s\) for every \(s>0\).
\end{proof}

\section{Prescribed discrepancy by inserting random points}

We prove Theorem~\ref{thm:prescribed-intermediate} by inserting
independent uniform points in \([0,1/2)\) into a sequence satisfying
Theorem~\ref{thm:low-discrepancy-background}.
Their number determines the discrepancy asymptotic, while estimates
for pairs involving the inserted points yield PPC almost surely.

\subsection{Construction of the sequence}

Let \(\psi\) satisfy the assumptions of
Theorem~\ref{thm:prescribed-intermediate}, and put \(h_N=N\psi(N)\),
with \(h_0=0\). Choose \(\kappa>0\) such that \(\kappa L\le1/2\), where $L$ is the increment bound in
\eqref{eq:intermediate-increments}. Set
\[
H_0=0,
\qquad
H_N=\lfloor\kappa h_N\rfloor,
\qquad
M_N=N-H_N
\quad(N\ge1).
\]
From \eqref{eq:intermediate-increments} we then obtain
\[
H_N-H_{N-1}\in\{0,1\},
\qquad
0\le H_N\le N/2
\qquad(N\ge1).
\]
Since \(H_N=\kappa h_N+O(1)\), we deduce from
\eqref{eq:intermediate-growth-lower} and
\eqref{eq:intermediate-growth-upper} that
\[
H_N\ll\sqrt{N\log N},
\qquad
\frac{H_N}{\log N}\to\infty,
\qquad
\frac{M_N}{N}\to1.
\]
Consider the probability space \(( \Omega, {\mathcal F}, {\PP})\), where
\[
  \Omega=[0,1/2)^{\mathbb N},
\qquad
 {\mathcal F}=\mathcal B([0,1/2))^{\otimes\mathbb N},
\qquad
 {\PP}=\bigl(2\lambda|_{[0,1/2)}\bigr)^{\otimes\mathbb N}.
\]
For \(\omega=(z_1,z_2,\ldots)\in \Omega\), we define
\[
W_n(\omega)=z_n
\qquad(n\ge1).
\]
Then \(W_1,W_2,\ldots\) are independent and uniform on \([0,1/2)\).
Fix a sequence \((y_m)_{m\ge1}\) satisfying the conclusions of
Theorem~\ref{thm:low-discrepancy-background}.
We insert \(W_1,W_2,\ldots\) into this sequence, retaining all its terms
in their original order. More precisely, we define
\begin{equation}
\notag
X_N=
\begin{cases}
W_{H_N},&H_N-H_{N-1}=1,\\
y_{M_N},&H_N-H_{N-1}=0.
\end{cases}
\end{equation}
The first \(N\) terms thus consist of \(y_1,\ldots,y_{M_N}\)
and \(W_1,\ldots,W_{H_N}\).

\subsection{Discrepancy}

The distribution function of \(W_1\) is
\[
F_W(t)=\min\{2t,1\}
\qquad(0\le t\le1).
\]
For \(0\le t\le1\), from the construction we see
\begin{align}
A_N^X(t)-Nt
&=\bigl(A_{M_N}^y(t)-M_Nt\bigr)
+\bigl(A_{H_N}^W(t)-H_NF_W(t)\bigr)
\notag\\
&\quad+H_N\min\{t,1-t\}.
\notag
\end{align}
The Dvoretzky--Kiefer--Wolfowitz inequality in Massart's sharp form
\cite{Massart1990} gives, for every integer \(H\ge2\),
\[
 {\PP}\left(
\sup_{0\le t\le1}|A_H^W(t)-HF_W(t)|>\sqrt{H\log H}
\right)\le\frac2{H^2}.
\]
By the Borel--Cantelli lemma, we have, almost surely,
\begin{equation}
\label{eq:intermediate-empirical-bound}
\sup_{0\le t\le1}|A_H^W(t)-HF_W(t)|
=O\bigl(\sqrt{H\log H}\bigr)
\qquad(H\to\infty).
\end{equation}
Since \(\sup_{0\le t\le1}\min\{t,1-t\}=1/2\), we deduce from the reverse triangle
inequality, \eqref{eq:main-low-discrepancy} and
\eqref{eq:intermediate-empirical-bound} that, almost surely,
\begin{align}
\left|ND_N^*-\frac{H_N}{2}\right|
&\le
\sup_{0\le t\le1}|A_{M_N}^y(t)-M_Nt|
+\sup_{0\le t\le1}|A_{H_N}^W(t)-H_NF_W(t)|
\notag\\
&\ll\log N+\sqrt{H_N\log H_N}.
\notag
\end{align}
Using \(H_N=\kappa h_N+O(1)\), we then obtain, almost surely,
\begin{equation}
\notag
ND_N^*=\frac\kappa2h_N
+O\left(\log N+\sqrt{h_N\log h_N}\right).
\end{equation}
By \eqref{eq:intermediate-growth-lower}, the error is \(o(h_N)\), so we conclude that
 \(
D_N^*\sim (\kappa/2) \psi(N).\)

\subsection{Poissonian pair correlation}

Fix \(s>0\). The following argument applies to every realization
satisfying \eqref{eq:intermediate-empirical-bound}.
Splitting the ordered pairs according to their origins, we obtain
\begin{equation}
\notag
\begin{aligned}
R_N^X(s)
={}&\frac{M_N}{N}
R_{M_N}^y\left(s\frac{M_N}{N}\right)+\frac2N\sum_{j=1}^{H_N}
\sharp\left\{1\le m\le M_N:
\normT{y_m-W_j}\le\frac{s}{N}\right\}\\
&+\frac1N
\sharp\left\{1\le i\ne j\le H_N:
\normT{W_i-W_j}\le\frac{s}{N}\right\}.
\end{aligned}
\end{equation}
The first term tends to \(2s\), by the PPC of \((y_m)\),
monotonicity in the window parameter, and \(M_N/N\to1\).
Applying \eqref{eq:main-low-discrepancy} and
\eqref{eq:intermediate-empirical-bound} to circular intervals of radius
\(s/N\), we obtain, uniformly in \(z\in\mathbb T\),
\[
\begin{aligned}
\sharp\left\{1\le m\le M_N:\normT{y_m-z}\le\frac{s}{N}\right\}
&\ll_s\log N,\\
\sharp\left\{1\le j\le H_N:\normT{W_j-z}\le\frac{s}{N}\right\}
&\ll_s\sqrt{H_N\log H_N}.
\end{aligned}
\]
Since \(H_N\ll\sqrt{N\log N}\), we deduce that
\[
\begin{aligned}
0
\le R_N^X(s)
-\frac{M_N}{N}R_{M_N}^y\left(s\frac{M_N}{N}\right)\ll_s
\frac{H_N\log N}{N}
+\frac{H_N^{3/2}\sqrt{\log H_N}}{N}
\to 0.
\end{aligned}
\]
Thus \((X_n)\) has PPC on the same probability-one event.

\section{Prescribed discrepancy by independent sampling}
\label{sec:prescribed-bias}

We prove Theorem~\ref{thm:prescribed-above-square-root} using independent
random variables with nonidentical distributions.
We control discrepancy through their average density and prove PPC
almost surely by estimating the expectation and variance of \(R_N(s)\).

\subsection{Construction of the sequence}

Let $\psi$ satisfy the assumptions of
Theorem~\ref{thm:prescribed-above-square-root}. Put
\(
 h_N=N\psi(N),\) with \(h_0=0.
\)
Choose $\kappa>0$ such that
\(
\kappa L\le\1/2
\).
 Define
\(
 a_n=\kappa(h_n-h_{n-1}).
\)
Then \(0\le a_n\le 1/2,\)
and
\begin{equation}
\label{eq:an-sum}
\sum_{n=1}^N a_n=\kappa h_N=\kappa N\psi(N).
\end{equation}

Define  
\(
g(x)=
\1_{[0,1/2)}(x)-\1_{[1/2,1)}(x)
\), and, for  $n\ge1$, define 
\(
 f_n(x)=1+a_ng(x).
\)
We see that 
\(
1/2\le f_n(x)\le 3/2,
\)
and since \(\int_{\mathbb T}g(x)\dd x=0\), each \(f_n\) is a
probability density.
Consider the probability space \(( \Omega, {\mathcal F}, \PP)\), where
\[
 \Omega=\mathbb T^{\mathbb N},
\qquad
 {\mathcal F}=\mathcal B(\mathbb T)^{\otimes\mathbb N},
\qquad
 \PP=\bigotimes_{n\ge1}(f_n(x)\dd x).
\]
For \(\omega=(x_1,x_2,\ldots)\in \Omega\), we define
\[
X_n(\omega)=x_n
\qquad(n\ge1).
\]
Then \(X_1,X_2,\ldots\) are independent random variables
with respective densities \(f_1,f_2,\ldots\).

Their average density is
\begin{equation}
\label{eq:average-biased-density}
\overline f_N(x)
=\frac1N\sum_{n=1}^N f_n(x)
=1+\kappa\psi(N)g(x).
\end{equation}

\subsection{Discrepancy}

For \(0\le t\le1\), write
\[
G_N(t)=\frac1N\sum_{n=1}^N\1_{[0,t)}(X_n),
\qquad
F_N(t)=\mathbb E G_N(t)=\int_0^t\overline f_N(x)\dd x.
\]
Since \(0\le\overline f_N\le3/2\), the monotonicity of \(G_N\) yields
\[
\sup_{0\le t\le1}|G_N(t)-F_N(t)|
\le
\max_{0\le j\le N}
|G_N(j/N)-F_N(j/N)|+\frac{3}{2N}.
\]
We combine Hoeffding's inequality \cite[Theorem~1]{Hoeffding1963}
with a union bound to obtain, for every \(\varepsilon>0\),
\[
\PP\left(
\max_{0\le j\le N}
|G_N(j/N)-F_N(j/N)|>\varepsilon
\right)
\le2(N+1)e^{-2N\varepsilon^2}.
\]
Taking \(\varepsilon=2\sqrt{(\log N)/N}\) and applying the
Borel--Cantelli lemma, we obtain that, almost surely,
\begin{equation}
\label{eq:empirical-almost-sure}
\sup_{0\le t\le1}|G_N(t)-F_N(t)|
=O\left(\sqrt{\frac{\log N}{N}}\right).
\end{equation}

By \eqref{eq:average-biased-density},
\[
F_N(t)-t
=\kappa\psi(N)\int_0^tg(x)\dd x
=\kappa\psi(N)\min\{t,1-t\}
\qquad(0\le t\le1).
\]
Hence
\begin{equation}
\notag
\sup_{0\le t\le1}|F_N(t)-t|=\frac\kappa2\psi(N).
\end{equation}
The reverse triangle inequality for the supremum norm yields
\[
\left|D_N^*-\frac\kappa2\psi(N)\right|
\le\sup_{0\le t\le1}|G_N(t)-F_N(t)|.
\]
Equation~\eqref{eq:empirical-almost-sure} therefore implies that, almost surely,
\begin{equation}
\notag
ND_N^*=\frac\kappa2h_N
+O\left(\sqrt{N\log N}\right).
\end{equation}
By \eqref{eq:above-empirical-scale}, the error is $o(h_N)$, so we conclude that
$D_N^*\sim(\kappa/2)\psi(N)$.

\subsection{Poissonian pair correlation}

Fix \(s>0\). For \(N>2s\) and \(1\le i\ne j\le N\), put
\[
B_{ij}=\1_{\{\normT{X_i-X_j}\le s/N\}}.
\]
Then
\[
R_N(s)=\frac1N\sum_{1\le i\ne j\le N}B_{ij}
=\frac2N\sum_{1\le i<j\le N}B_{ij}.
\]
By the independence, expanding \(f_i(x)f_j(y)\) and using
\(\int_{\mathbb T}g(x)\dd x=0\), we obtain
\[
\mathbb E B_{ij}-\frac{2s}{N}
=
a_i a_j\int_{\mathbb T}
\int_{\normT{x-y}\le s/N}g(x)g(y)\dd y\dd x.
\]
Since \(|g|\le1\) and \(a_i,a_j\ge0\),
\[
\left|\mathbb E B_{ij}-\frac{2s}{N}\right|
\le\frac{2s}{N}a_i a_j.
\]
Summing over \(i\ne j\) and using \eqref{eq:an-sum} we obtain
\begin{align}
\left|\mathbb E(R_N(s)-2s)\right|
 \le
\frac{2s}{N}
+\frac{2s}{N^2}\sum_{1\le i\ne j\le N}a_i a_j
\le
\frac{2s}{N}+2s\kappa^2\psi(N)^2
\to 0.
\notag
\end{align}

For pairwise distinct indices \(i,j,k\), conditioning on \(X_i\),
using the independence and \(f_n\le3/2\), and then taking expectations, we have 
\[
\mathbb E B_{ij}\le\frac{3s}{N},
\qquad
\mathbb E(B_{ij}B_{ik})\le\frac{9s^2}{N^2}.
\]
Indicators corresponding to disjoint index sets are independent.
Using \(\operatorname{Var}(B_{ij})\le\mathbb E B_{ij}\) and
\(\operatorname{Cov}(B_{ij},B_{ik})\le\mathbb E(B_{ij}B_{ik})\),
we obtain
\begin{align}
\operatorname{Var}(R_N(s))
&=\frac4{N^2}
\operatorname{Var}\!\left(\sum_{1\le i<j\le N}B_{ij}\right)
\notag\\
&\le\frac4{N^2}
\left(
\binom N2\frac{3s}{N}
+N(N-1)(N-2)\frac{9s^2}{N^2}
\right)
\le\frac{6s+36s^2}{N}.
\notag
\end{align}

The standard subsequence argument now shows that
\((X_n)\) has PPC almost surely.

\section{Optimal discrepancy constants}
\label{sec:consequences-scope}

Theorem~\ref{thm:low-discrepancy-background} shows that PPC is compatible
with the optimal discrepancy order $(\log N)/N$.
We now consider the optimal limsup constant. Set
\[
c^*:=\inf_x\limsup_{N\to\infty}
\frac{ND_N^*(x)}{\log N},
\]
where the infimum is over all infinite sequences $x=(x_n)_{n\ge1}$
in $[0,1)$. Define
\[
c^*_{\mathrm{PPC}}:=
\inf_{\substack{x\\\text{$x$ has PPC}}}
\limsup_{N\to\infty}\frac{ND_N^*(x)}{\log N}.
\]
Since PPC sequences form a subclass of all sequences,
\begin{equation}
\notag
c^*\le c^*_{\mathrm{PPC}}.
\end{equation}
The exact value of $c^*$ is unknown; for bounds,  we refer to
\citet{LarcherPuchhammer2016} and \citet{Ostromoukhov2009}.

\begin{openproblem}
Is $c^*_{\mathrm{PPC}}=c^*$?
\end{openproblem}

\section*{Statements and Declarations}

\noindent\textbf{Funds.}
This work is supported by National Key R\&D Program of China
(No. 2024YFA1013700) and NSFC (Nos. 12671113 and 12331005),
and the Fundamental Research Funds for the Central Universities.

\medskip
\noindent\textbf{Competing Interests.}
The authors declare that they have no financial or non-financial
interests relevant to the content of this article.

\medskip
\noindent\textbf{Data Availability.}
This article is theoretical. No research datasets were generated
or analysed.

\bigskip

\noindent School of Mathematics and Statistics, Huazhong University of
Science and Technology, Wuhan 430074, China.

\noindent\textit{Email address:} \texttt{tanbo@hust.edu.cn}

\bigskip

\noindent School of Mathematics and Statistics, Huazhong University of
Science and Technology, Wuhan 430074, China.

\noindent\textit{Email address:} \texttt{zhangzhen\_@hust.edu.cn}

\bigskip

\noindent School of Mathematics and Statistics, Wuhan University of
Technology, Wuhan 430070, China.

\noindent\textit{Email address:} \texttt{zhouql@whut.edu.cn}

\end{document}